\documentclass[a4paper]{amsart}
\usepackage{latexsym}
\usepackage{amsmath}
\usepackage[english]{babel}
\usepackage{amssymb}

\newtheorem{theorem}{Theorem}[section]

\newtheorem{lemma}[theorem]{Lemma}

\newtheorem{definition}[theorem]{Definition}
\newtheorem{definition and theorem}[theorem]{Definition and Theorem}
\newtheorem{remark}[theorem]{Remark}
\newtheorem{*remark}[theorem]{$^* $Remark}

\newtheorem{*exercise}[theorem]{$^* $Exercise}
\newtheorem{**exercise}[theorem]{$^{** } $Exercise}

\newtheorem{corollary}[theorem]{Corollary}

\begin{document}
\title[VaR and Expected Shortfall for Mixture Elliptic Linear Portfolios]{VaR and ES For Linear Portfolios with Mixture of Elliptic Distributed Risk Factors}
\author[Jules SADEFO KAMDEM]{Jules SADEFO KAMDEM \\ \newline Universit\'e de Reims \\ Universit\'e d'Evry\\}
\thanks{This draft is a part  of J.SADEFO-KAMDEM PhD Thesis at the Universit\'e de Reims (France).
The author is a temporary lecturer at the Department of
Mathematics ( Universit\'e d'evry Val d'essonne France).
\underline{Author's address}: Universit\'e de Reims, Laboratoire
de Math\'ematiques UMR 6056-CNRS , BP 1039 Moulin de la
Housse , 51687 Reims cedex 2 FRANCE.\\
e-mail: sadefo@univ-reims.fr}

\maketitle

\begin{abstract}
In this paper, following the generalization of {\em Delta Normal VaR} to {\em Delta Mixture Elliptic VaR}
in  Sadefo-Kamdem \cite{S2}, we give and explicit formula to estimate linear VaR and
ES when the risk factors changes with the mixture of $\it{t}$-Student
distributions. In particular, we give rise to {\em Delta-Mixture-Student VaR} and the {\em Delta-Mixture-Elliptic ES}.
\end{abstract}

\bigskip


\section{{\bf Introduction} }

 \medskip

The original RiskMetrics methodology for estimating VaR was based
on parametric methods, and used the multi-variate normal
distribution. This approach works well for the so-called {\em
linear portfolios}, that is, those portfolios whose aggregate
return is, to a good approximation, a linear function of the
returns of the individual assets which make up the portfolio, and
in situations where the latter can be assumed to be jointly
normally distributed. For other portfolios, like portfolios of
derivatives depending non-linearly on the return of the
underlying, or portfolios of non-normally distributed assets, one
generally turns to Monte Carlo methods to estimate the VaR. This
is an issue in situations demanding for real-time evaluation of
financial risk. For non-linear portfolios, practitioners, as an
alternative to Monte Carlo, use $\Delta $-normal VaR methodology,
in which the portfolio return is linearly approximated, and an
assumption of normality is made. Such methods present us with a
trade-off between accuracy and speed, in the sense that they are
much faster than Monte Carlo, but are much less accurate unless
the linear approximation is quite good and the normality
hypothesis holds well. The assumption of normality simplifies the
computation of VaR considerably. However it is inconsistent with
the empirical evidence of assets returns, which finds that asset
returns are fat tailed. This implies that extreme events are much
more likely to occur in practice than would be predicted based on
the assumption of normality.

Some alternative return distributions have been proposed in
the world of elliptic distributions by Sadefo-Kamdem
\cite{S2}, that better reflect the empirical evidence. In
this paper, following \cite{S2}, I examine one such
alternative that simultaneously allows for asset returns
that are fat tailed and for tractable calculation of Value-at-Risk
and Expected Shortfall, by giving attention to
mixture of elliptic distributions, with an explicit formula  of VaR and ES in the
special case where assets returns changes with mixture of Student-$\it{t}$
distributions. Note that, the particular case based on  mixture of normal
distributions, has been proposed by Zangari(1996)\cite{Z}, Subu-Venkataraman \cite{SV}
and some references therein.
\medskip

An obvious first generalization is to keep the linearity
assumption, but replace the normal distribution by some other
family of multi-variate distributions. In Sadefo-Kamdem \cite{S2},
we have such generalization concerning linear portfolios, in the
case where the joint risk factors changes with mixture of elliptic
distributions. In this paper, by using some generalized theorems
concerning {\em Delta-Mixture-Elliptic} VaR and {\em
Delta-Student} VaR  in \cite{S2}, we introduce the notion of {\em
Delta-Mixture-Student VaR}, {\em Delta-Mixture-Elliptic ES} and
the {\em Delta-Mixture-Student ES}.

So the particular subject of this paper, is to give an explicit
formulas that will permit to obtain the linear VaR or linear ES,
when the joint risk factors of the linear portfolios, changes with
mixture of $\it{t}$-Student distributions.   Note that, since one
shortcoming of the multivariate t-distribution is that all the
marginal distributions must have the same degrees of freedom,
which implies that all risk factors have equally heavy tails, the
mixture of $\it{t}$-Student will be view as a serious
alternatives, to a simple $\it{t}$-Student-distribution.
Therefore, the methodology proposes by this paper seem to be
interesting to controlled thicker tails than the standard Student
distribution.

\medskip \quad The paper
is organized, as follows: In section 2, we recall some theorems
concerning the {\em Delta-Elliptic}, {\em Delta-Elliptic} and {\em
Delta Mixture Elliptic VaR} given by  Sadefo-Kamdem \cite{S2}. In
section 3,  following the theorem concerning {\em Delta Mixture
Elliptic VaR}, we show how to reduce the computation of the {\em
Delta-Mixture-Student VaR} to finding the zeros of a mixture of
special function. In section 4, we introduces the notion of {\em
Delta Mixture Elliptic ES },  by treat the expected shortfall for
general mixture of elliptic distribution,  with special attention
to {\em Delta Mixture Elliptic ES }. Finally, in section 5 we
discuss some potential application areas.

\section{Some Notions on log-elliptic Linear VaR}
\medskip

In this section, following \cite{S2}, we recall some notions on
elliptic distributions and Linear VaR.
\medskip

We will use the following notational conventions for the action of
matrices on vectors: single letters $x , y , \cdots $ will denote
{\em row vectors} $(x_1 , \cdots , x_n ) $, $(y_1 , \cdots y_n )
$. The corresponding column vectors will be denoted by $x^t , y^t
$,the $^t $ standing more generally for taking the transpose of
any matrix. Matrices $A = (A_{ij } )_{i, j } $, $B $ , etc. will
be multiplied in the usual way. In particular, $A $ will act on
vectors by left-multiplication on column vectors, $A y^t $, and by
right multiplication on row vectors, $xA $; $x \cdot x = x x^t =
x_1 ^2 + \cdots + x_n ^2 $ will stand for the Euclidean inner
product.
\medskip

A portfolio with time-$t $ value $\Pi (t) $ is called linear if
its profit and loss $\Delta \Pi(t)= \Pi (t) - \Pi (0) $ over a
time window, [0 t] is a linear function of the returns
$X_{1}(t),\ldots,X_{n}(t)$ of its constituents over the same time
period:
$$
\Delta \Pi(t)= \delta_{1} X_{1}+\delta_{2} X_{2}+...+\delta_{n} X_{n}
$$
This will for instance be the case for ordinary portfolios of
common stock, if we use percentage returns, and will also hold to
good approximation with log-returns, provided the time window
[0,t] is small. We will drop the time $t $ from our notations,
since it will be kept fixed, and simply write $X_{j}$,$\Delta
\Pi$, etc. We also put
$$
X = (X_1 , \cdots , X_n ) ,
$$
so that $\Delta \Pi = \delta \cdot X = \delta X^t . $
\medskip

We now assume that the $X_j $ are elliptically distributed with
mean $\mu$ and correlation matrix $\Sigma=AA^t$:
$$
(X_{1},\ldots,X_{n})\sim N(\mu,\Sigma , \phi ).
$$

This means that the pdf of X is of the form
$$
f_{X}(x)=|\Sigma|^{-1/2} g((x-\mu)\Sigma^{-1}(x-\mu)^t ) ,
$$
where $|\Sigma | $ stands for the determinant of $\Sigma $, and
where $g: \mathbb{R }_{\geq 0 } \to 0 $ is such that the Fourier
transform of $g(|x|^2 ) $, as a generalized function on $\mathbb{R}^n $,
 is equal to $\phi (|\xi |^2 ) $

\footnote{One uses $\phi $
as a parameter for the class of elliptic distributions, since it
is always well-defined as a continuous function: $\phi (|\xi |^2 )$
 is simply the characteristic function of an $X \sim N(0, Id ,\phi ) $.
Note, however, that in applications we'd rather know $g$}.

Assuming that $g $ is continuous, and non-zero everywhere, the
Value at Risk at a confidence level of $1 - \alpha$ is given by
solution of the following equation:

Here we follow
the usual convention of recording portfolio losses by negative
numbers, but stating the Value-at-Risk as a positive quantity of
money.

\subsection{Linear VaR with mixtures of elliptic Distributions}
\medskip
\quad
       Mixture distributions can be used to model situations where
the data can be viewed as arising from two or more distinct
classes of populations; see also \cite{MX}. For example, in the
context of Risk Management, if we divide trading days into two
sets, quiet days and hectic days, a mixture model will be based on
the fact that returns are moderate on quiet days, but can be
unusually large or small on hectic days. Practical applications of
mixture models to compute VaR can be found in Zangari \cite{Z}
(1996), who uses a mixture normal to incorporate fat tails in VaR
estimation. In Sadefo-Kamdem \cite{S2}, we have generalized the
preceding section to the situation where the joint log-returns
follow a mixture of elliptic distributions, that is, a convex
linear combination of elliptic distributions. In this section, a
special attention will be give to mixture of Student-$\it{t}$
distributions.

\medskip
\begin{definition}\label{VaR-elliptic} \rm{We say that $ (X_{1},...,X_{n})$ has a joint distribution
that is the mixture  of $m $ elliptic distributions
$N(\mu_{j},\Sigma_{j},\phi_{j}) $\footnote{or
$N(\mu_{j},\Sigma_{j},g_{j})$ if we parameterize elliptical
distributions using $g $ instead of $\phi $}, with weights
$\{\beta_{j}\}$  (j=1,..,m ;  $\beta_{j} > 0$ ;  $\sum_{j=1}^m
\beta_{j} = 1$), if its cumulative distribution function can be
written as
$$
F_{X_{1},...,X_{n}}(x_{1},...,x_{n}) = \sum_{j=1}^m \beta_{j}
F_{j}(x_{1},...,x_{n})
$$
with $F_{j}(x_{1},...,x_{n}) $ the cdf of
$N(\mu_{j},\Sigma_{j},\phi_{j}) $. } \end{definition}

\begin{remark} \rm{
In practice, one would usually limit oneself to $m = 2 $, due to
estimation and identification problems; see \cite{MX}. }
\end{remark}

The following lemma is given by Sadefo-Kamdem \cite{S2}.
\medskip
\begin{lemma}
Let $\Delta \Pi=\delta_{1} X_{1}+\ldots+\delta_{n} X_{n}$  with
 $(X_{1},\ldots,X_{n})$ is a mixture of elliptic distributions, with
 density
$$
f(x)=\sum_{j=1}^m  \beta _j {|\Sigma_{j}|}^{-1/2 }
g_{j}((x-\mu_{j} )\Sigma_{j}^{-1}(x-\mu_{j} )^{t}) \label{dens}
$$
where $\mu_{j}$ is the vector mean, and  $\Sigma_{j} $ the
variance-covariance matrix of the $j $-th component of the
mixture.  We suppose that each $g_j $ is  integrable function over
$\mathbb{R}$, and that the $g_j $ never vanish jointly in a point
of $\mathbb{R }^m $. Then the value-at-Risk, or {\em Delta
mixture-elliptic VaR}, at confidence $1 - \alpha $ is given as the
solution of the transcendental equation
\begin{equation}
\alpha = \sum _{j = 1 } ^m \beta _j G_j \left( \frac{\delta
.\mu_{j}^{t} + VaR_{\alpha }}{(\delta \Sigma _j \delta )^{1/2 } }
\right) ,\label{eq3}
\end{equation}
where

$$G_j  = \frac{|S_{n-2}|}{
{|\Sigma_{j}|}^{1/2}} \int_{0}^{+\infty} r^{n-2} \Big{[}
\int_{-\infty }^\frac{- \delta \cdot \mu _j -VaR_{\alpha
}}{|\delta A_{j}| }      g_{j}( z_{1}^2 + r^{2}  ) dz_{1}\Big{]}
dr. \label{mixture}$$
  Here $ \delta
=(\delta_{1},\ldots,\delta_{n})$.
\end{lemma}
\begin{remark} \rm{
In the case of a mixture of m elliptic distributions the VaR will
not depend any more in a simple way on the total portfolio mean
and variance-covariance. This is unfortunate, but already the case
for a mixture of normal distributions. } \end{remark}
\begin{remark} \rm{
One might, in certain situations, try to model with a mixture of
elliptic distributions which all have the same variance-covariance
and the same mean, and obtain for example a mixture of different
tail behaviors by playing with the $g_j $'s. In that case the VaR
again simplifies to: $ VaR _{\alpha } = - \delta \cdot \mu +
q_{\alpha } \cdot \sqrt{ \delta \Sigma \delta ^t } $, with
$q_{\alpha } $ now the unique positive solution to
$$
\alpha  = \sum_{j=1}^m \beta_{j} G_j (q_{\alpha} ) .
$$
}
\end{remark}

\noindent The preceding can immediately be specialized to a
mixture of normal distributions. the details is left to the
reader.

\section{VaR with mixture  Student-$t$ distributions}

We now consider in detail the case where our mixture of elliptic
distributions is a mixture of multivariate  Student-$t $. We will,
unsurprisingly, call the corresponding $VaR$ the {\em Delta
mixture-Student VaR}.

In the case of our mixture of multi-variate t-Student
distributions, the portfolio probability density function is given
by: 
\begin{equation} 
h_X (x) =\sum_{j=1}^m \beta_{j}  \frac{\Gamma (\frac{\nu_{j} + n}{2})}{\Gamma(\nu_{j}/2).\sqrt{|\Sigma_{j}|(\nu_{j} \pi)^n }}
{\Big{(}1+\frac{(x-\mu_{j})^{t}\Sigma_{j}^{-1}(x-\mu_{j})}{\nu_{j}} \Big{)}}^{(\frac{-\nu_{j}-n}{2})} ,  \label{density}
\end{equation}
$x \in \mathbb{R}^{n} $ and $ \nu_{j} > 2   $. Hence $g_{j} $ is
given by
$$g_{j}(s)= C(\nu_{j} ,n) {(1+s/\nu_{j} )}^{-\frac{(n+\nu_{j} )}{2}} , \ \ s \geq 0 , $$
where we have put
$$
C(\nu_{j} ,n)=\frac{\Gamma (\frac{\nu_{j} +
n}{2})}{\Gamma(\nu_{j}/2) \sqrt{(\nu_{j} \pi)^n }} .
$$
Using this $g_{j} $ in $(\ref{eq3})$, we find that

\begin{equation} G_{j} (s) = \frac{\nu_{j}^\frac{n+\nu_{j}}{2} }{2 } |S_{n - 2
} | C(\nu_{j} , n )\int_s ^{\infty} I_{j}(z_1 ) dz_1 , \label{4j}
\end{equation}

where we have put \begin{equation} I_{j}(z_{1})=  \int_{z_{1}^2}^{+\infty}
(u-z_{1}^2)^{\frac{n-3}{2}} (\nu_{j} + u)^{-\frac{(n+\nu_{j}
)}{2}} du . \label{5j} \end{equation}

Following \cite{S2}, we have the following expression
\begin{eqnarray}
 G_{j}(s) = \frac{1 }{\nu_{j} \sqrt{\pi } } \left( \frac{\nu_{j} }{s^2 } \right)
 ^{\nu_{j} / 2 } \frac{\Gamma \left( \frac{\nu_{j} + 1 }{2 } \right) }
 {\Gamma \left( \frac{\nu_{j} }{2 } \right) }
{_2F}_{1}\Big{(}\frac{1+\nu_{j}}{2},\frac{\nu_{j}}{2};1+\frac{\nu_{j}}{2};-\frac{\nu_{j}
}{s^2 } \Big{)}  \label{9aj}
\end{eqnarray}
 we obtain the following corollary
\begin{corollary}\label{theomixtstud}
Let $\Delta \Pi=\delta_{1} X_{1}+\ldots+\delta_{n} X_{n}$  with
 $(X_{1},\ldots,X_{n})$ is a mixture of m $t$-Student distributions, with
 density $h_X$ defined by (\ref{density}),where $\mu_{j}$ is the vector mean,
  and  $\Sigma_{j} $ the variance-covariance matrix of the $j $-th component of the
mixture.   Then the value-at-Risk, or {\em Delta mixture-student
VaR}, at confidence $1 - \alpha $ is given as the solution of the
transcendental equation
\begin{equation}
\alpha  =\sum _{j = 1 } ^m  \frac{\beta _j \Gamma \left(
\frac{\nu_{j} + 1 }{2 } \right) }{\nu_{j} \sqrt{\pi }\Gamma \left(
\frac{\nu_{j} }{2 } \right) } {\nu_{j}}^{\frac{\nu_{j}}{2}} \left(
\frac{\delta .\mu_{j}^{t} + VaR_{\alpha } }{\delta \Sigma_{j}
\delta } \right) ^{\frac{-\nu_{j}}{2}}
{_2F}_{1}\Big{(}\frac{1+\nu_{j}}{2},\frac{\nu_{j}}{2};1+\frac{\nu_{j}}{2};-\frac{\nu_{j}(\delta
.\mu_{j}^{t} + VaR_{\alpha }) }{\delta \Sigma _j \delta } \Big{)}
\end{equation}
where $G_j $ is defined by (\ref{4j}) with $g = g_j . $ Here
$\delta =(\delta_{1},\ldots,\delta_{n})$.
\end{corollary}

\begin{corollary} \label{CorrmixtVaR2}
One might, in certain situations, try to model with a mixture of
$t$-Student distributions which all have the same
variance-covariance $\Sigma=\Sigma_{j}$ and the same mean
$\mu=\mu_{j}$, and obtain for example a mixture of different tail
behaviors by playing with the $\nu_{j} $'s. In that case the VaR
again simplifies to: $$ VaR _{\alpha } = - \delta \cdot \mu +
q_{\alpha } \cdot \sqrt{ \delta \Sigma \delta ^t }, $$ with
$q_{\alpha } $ now the unique positive solution to
$$
\alpha  =\sum _{j = 1 } ^m  \frac{\beta _j \Gamma \left(
\frac{\nu_{j} + 1 }{2 } \right) }{\nu_{j} \sqrt{\pi }\Gamma \left(
\frac{\nu_{j} }{2 } \right) } {\nu_{j}}^{\frac{\nu_{j}}{2}} \left(
\frac{\delta .\mu^{t} + VaR_{\alpha } }{\delta \Sigma \delta }
\right)
 ^{\frac{-\nu_{j}}{2}}
{_2F}_{1}\Big{(}\frac{1+\nu_{j}}{2},\frac{\nu_{j}}{2};1+\frac{\nu_{j}}{2};-\frac{\nu_{j}(\delta
.\mu^{t} + VaR_{\alpha }) }{\delta \Sigma \delta } \Big{)} .
$$
\end{corollary}

\begin{remark} \rm{
One might, in certain situations, try to model with a mixture of
$t$-Student distributions which all have the same $\nu_{j}=\nu$
and the same mean $\mu_{j}\approx 0$, and obtain for example a
mixture of different tail behaviors by playing with the
$\Sigma_{j} $'s. In that case the VaR is the unique positive
solution to
$$
\alpha  =\frac{\Gamma \left( \frac{\nu + 1 }{2 } \right) }{\nu
\sqrt{\pi }\Gamma \left( \frac{\nu }{2 } \right) }\sum _{j = 1 }
^m  \beta _j  \left( \frac{\nu ( VaR_{\alpha }) }{\delta \Sigma _j
\delta } \right)
 ^{\frac{\nu}{2}}
{_2F}_{1}\Big{(}\frac{1+\nu}{2},\frac{\nu}{2};1+\frac{\nu}{2};-\frac{\nu(
VaR_{\alpha }) }{\delta \Sigma _j \delta } \Big{)} .
$$
}

\end{remark}

\subsection{Some Numerical Result of Delta Mixture-Student VaR coefficient}
\medskip

Here we give some numerical results when applying the corollary
\ref{CorrmixtVaR2}, in the situation where $m=2$.

By introducing the function $F$ such that
 \begin{equation}
F(s,\beta,\nu_{1},\nu_{2})= \beta \cdot  G_{1}(s) + (1-\beta)\cdot
G_{2}(s), \end{equation} where $G_{j}$ is define in (\ref{9aj}), for $j=1,2$,
for given  as inputs $\beta$, $\nu_{1}$ and $\nu_{2}$, we  give a
table that contains some solutions
$s=q_{\beta,\nu_{1},\nu_{2}}=q_{\alpha}^{MS-VaR}$ of the following
transcendental equation:
$$F(s,\beta,\nu_{1},\nu_{2})=\alpha.\label{mixtwostu}$$
For given $\Sigma$, $\mu$, and $\delta$, these solutions will
permit to calculate $VaR_{\alpha}$, when the confidence is
$1-\alpha$.

\newpage
\begin{enumerate}
\item In the case where $\alpha=0.01$, we obtain some solutions of
(\ref{mixtwostu}) in the following table:
\begin{center}
\begin{tabular}{|c|c|c|c|c|c|c|c|c|}

  \hline

$(\nu_{1},\nu_{2})$  & (2,3)  & (3,4) & (4,6) & (5,8) & (6,10) &
(7, 15) & (8, 40) & (9, 16)
  \\
\hline
  $q_{0.05,\nu_{1},\nu_{2}}$   & 4.64839 & 3.78507  & 3.17184 & 3.91919 & 2.78228 & 2.62175  & 2.44602 & 2.59524  \\
 \hline
 $q_{0.10,\nu_{1},\nu_{2}}$  &4.7586 & 3.82348  & 3.20124 & 2.94213 & 2.80092 & 2.64116  & 2.46906 & 2.60704    \\
 \hline
 $q_{0.15,\nu_{1},\nu_{2}}$   & 4.87115 & 3.86216  & 3.23086 & 2.9652 & 2.81965 & 2.6607  & 2.49235 & 2.61887   \\
 \hline
 $q_{0.20,\nu_{1},\nu_{2}}$   & 4.98587 & 3.9011  & 3.26066 & 2.98846 & 2.83846 & 2.68035  & 2.51586 & 2.63073   \\
 \hline
 $q_{0.25,\nu_{1},\nu_{2}}$   & 5.10258 & 3.94025  & 3.29063 & 3.01177 & 2.85734 & 2.70009  & 2.53957 & 2.64261   \\
 \hline
 $q_{0.30,\nu_{1},\nu_{2}}$   & 5.22106 & 3.97962  & 3.32075 & 3.03518 & 2.87629 & 2.71991  & 2.56344 & 2.65452   \\
 \hline
$q_{0.35,\nu_{1},\nu_{2}}$   & 5.34113 & 4.01917  & 3.35100 & 3.05866 & 2.89528 & 2.7398  & 2.58744 & 2.66644   \\
 \hline
 $q_{0.40,\nu_{1},\nu_{2}}$   & 5.46259 & 4.05888  & 3.38136 & 3.08221 & 2.91432 & 2.75974  & 2.6115 & 2.67838   \\
 \hline
 $q_{0.45,\nu_{1},\nu_{2}}$   & 5.58523 & 4.09873  & 3.41180 & 3.10502 & 2.93339 & 2.77972  & 2.6357 & 2.69033   \\
 \hline
 $q_{0.50,\nu_{1},\nu_{2}}$   & 5.70886 & 4.13870  & 3.44231 & 3.12946 & 2.95248 & 2.79972  & 2.65989 & 2.70228   \\
 \hline
\end{tabular}
\\
\begin{tabular}{|c|c|c|c|c|c|c|c|c|}
  \hline

      $(\nu_{1},\nu_{2})$  & (10,20) & (20,30) & (200, 300) & (250,50) & (275,15) & (300,55) & (400,10) & (1000,5) \\
\hline
$q_{0.05,\nu_{1},\nu_{2}}$  & 2.53963 & 2.46079 & 2.33916 & 2.40018  & 2.58957 &  2.39322 & 2.7432 & 3.3202  \\
 \hline
  $q_{0.10,\nu_{1},\nu_{2}}$  & 2.55132  &2.46432 & 2.33947  & 2.39709 & 2.57661 & 2.39036 & 2.72242 & 3.27401   \\
 \hline
$q_{0.15,\nu_{1},\nu_{2}}$  & 2.56304  & 2.46785 & 2.33978 &  2.39399  & 2.56359 & 2.38750 & 2.7014 & 3.22632  \\
 \hline
 $q_{0.20,\nu_{1},\nu_{2}}$   & 2.5748 & 2.47139  & 2.3401 & 2.3909 & 2.55051 & 2.38464  & 2.68019 & 3.17715   \\
 \hline
 $q_{0.25,\nu_{1},\nu_{2}}$   & 2.58658 & 2.47492  & 2.34041 & 2.3878 & 2.53738 & 2.38178  & 2.6588 & 3.12651   \\
 \hline
 $q_{0.30,\nu_{1},\nu_{2}}$   & 2.59838 & 2.47846  & 2.34073 & 2.38471 & 2.52422 & 2.37892  & 2.63726 & 3.07446   \\
 \hline
 $q_{0.35,\nu_{1},\nu_{2}}$   & 2.6102 & 2.482  & 2.34104 & 2.38161 & 2.51102 & 2.37605  & 2.61559 & 3.02112   \\
 \hline
 $q_{0.40,\nu_{1},\nu_{2}}$   & 2.62204 & 2.48553  & 2.34136 & 2.37851 & 2.49779 & 2.37319  & 2.59382 & 2.96663   \\
 \hline
 $q_{0.45,\nu_{1},\nu_{2}}$   & 2.63389 & 2.48907  & 2.34167 & 2.37541 & 2.48455 & 2.37033  & 2.57198 & 2.91121   \\
 \hline
 $q_{0.50,\nu_{1},\nu_{2}}$   & 2.64574 & 2.49261  & 2.34199 & 2.37232 & 2.4713 & 2.36746  & 2.55009 & 2.85513   \\
 \hline
\end{tabular}
\end{center}
\medskip

\item In the case where $\alpha=0.001$, we obtain some solutions
of (\ref{mixtwostu}) in the following table:
\begin{center}
\begin{tabular}{|c|c|c|c|c|c|c|c|c|}

  \hline

$(\nu_{1},\nu_{2})$  & (2,3)  & (3,4) & (4,6) & (5,8) & (6,10) &
(7, 15) & (8, 40) & (9, 16)
  \\
\hline
   $q_{0.20,\nu_{1},\nu_{2}}$   & 12.8878 & 7.84891  & 5.66393 & 4.82769 & 4.39245 & 3.98902  & 3.62286 & 3.82625   \\
 \hline
 $q_{0.25,\nu_{1},\nu_{2}}$   & 13.5577 & 8.01412  & 5.77451 & 4.90665 & 4.45334 & 4.05064  & 3.69896 & 3.86013   \\
 \hline
 $q_{0.30,\nu_{1},\nu_{2}}$   & 14.2205 & 8.17734  & 5.88317 & 4.98414 & 4.51241 & 4.11084  & 3.77242 & 3.89346   \\
 \hline
$q_{0.35,\nu_{1},\nu_{2}}$   & 14.874 & 8.33840  & 5.98975 & 5.06004 & 4.57030 & 4.16948  & 3.84285 & 3.92621   \\
 \hline
 $q_{0.40,\nu_{1},\nu_{2}}$   & 15.5168 & 8.49717  & 6.09412 & 5.13427 & 4.62694 & 4.22648  & 3.91007 & 3.95838   \\
 \hline
 $q_{0.45,\nu_{1},\nu_{2}}$   & 16.1480 & 8.65357  & 6.19624 & 5.20677 & 4.68229 & 4.28179  & 3.97400 & 3.98993   \\
 \hline
 $q_{0.50,\nu_{1},\nu_{2}}$   & 16.7671 & 8.80753  & 6.29604 & 5.27752 & 4.73634 & 4.33537  & 3.03470 & 4.02087   \\
 \hline
\end{tabular}
\end{center}

\medskip

\item In the case where $\alpha=0$, we obtain some solutions of
(\ref{mixtwostu}) in the following table:
\begin{center}

\begin{tabular}{|c|c|c|c|c|c|c|c|c|}

  \hline

$(\nu_{1},\nu_{2})$  & (2,3)  & (3,4) & (4,6) & (5,8) & (6,10) &
(7, 15) & (8, 40) & (9, 16)
  \\
\hline
   $q_{0.20,\nu_{1},\nu_{2}}$   & 322.785 & 82.6688  & 31.0894 & 20.7154 & 15.8813 & 11.4371  & 10.1089 & 9.25604   \\
 \hline
 $q_{0.25,\nu_{1},\nu_{2}}$   & 352.09 & 87.1881  & 32.5561 & 21.541 &  16.42471  & 11.7949 & 10.3957 & 9.47529   \\
 \hline
 $q_{0.30,\nu_{1},\nu_{2}}$   & 378.302 & 91.2285  & 33.8309 & 22.2487 &  16.88721  &12.0958 &10.6352 & 9.66243   \\
 \hline
$q_{0.35,\nu_{1},\nu_{2}}$   & 402.155 & 94.8927  & 34.9619 & 22.8697 &  17.2907 &   12.3561 & 10.8414 & 9.82571 \\
 \hline
 $q_{0.40,\nu_{1},\nu_{2}}$   & 424.137 & 98.2529  & 35.981 & 23.4244 &  17.6493  & 12.5858 & 11.0227 & 9.97061  \\
 \hline
 $q_{0.45,\nu_{1},\nu_{2}}$   & 444.591 & 101.362  & 36.9102 & 23.9265 &  17.9726  & 12.7919 &11.1848 & 10.1009   \\
 \hline
 $q_{0.50,\nu_{1},\nu_{2}}$   & 463.771 & 104.26  & 37.7655 & 24.3858 &  18.2673  & 12.9789 & 11.3316 & 10.2194  \\
 \hline
\end{tabular}
\end{center}

\end{enumerate}

\begin{remark} \rm{ Note that, the precedent results are available when $\alpha=0$. This means that with our model,
one would calculate the linear VaR with mixture of elliptic
distributions, for 100 percent confidence level. }
\end{remark}

\section{Expected Shortfall with mixture of elliptic distributions}
\medskip

Expected shortfall is a sub-additive risk statistic that describes
how large losses are on average when they exceed the VaR level.
Expected shortfall will therefore give an indication of the size
of extreme losses when the VaR threshold is breached. We will
evaluate the expected shortfall for a linear portfolio under the
hypothesis of mixture of elliptically distributed risk factors.
Mathematically, the expected shortfall associated with a given VaR
is defined as:
$$
\mbox{Expected Shortfall } = \mathbb{E } (-\Delta \Pi \vert
-\Delta \Pi > VaR ),
$$
see for example \cite{MX}. Assuming again a multivariate mixture
of elliptic probability density $f_{X}(x) = \sum_{i=1}^n \beta_{i}
{|\Sigma_{i}|}^{-1/2} g_{i}((x-\mu_{i} )\Sigma_{i}^{-1}(x-\mu_{i}
)^{t}) $, the Expected Shortfall at confidence level $1 - \alpha $
is given by
\begin{eqnarray*}
- ES_{\alpha } &=& \mathbb{E } ( \Delta \Pi \mid   \Delta \Pi\leq
-VaR_{\alpha } ) \\
&=& \frac{1 }{\alpha } \mathbb{E } \left( \Delta \Pi \cdot 1_{\{
\Delta
\Pi \leq -VaR_{\alpha } \} } \right) \\
&=& \frac{1}{\alpha }\int_{\{ \delta x^t \leq -VaR_{\alpha }\}}
\delta x^t \ f_{X}(x) \ dx  \\
&=& \sum_{i=1}^{n}\beta_{i} \frac{ {|\Sigma_{i}|}^{-1/2}}{\alpha}
\int_{\{ \delta x^t \leq -VaR_{\alpha }\}} \delta x^t \ g_{i}
((x-\mu_{i} )\Sigma_{i}^{-1}(x-\mu_{i} )^{t}) dx .
\end{eqnarray*}
Let $\Sigma = A_{i}^t \; A_{i} $, as before.Doing the same linear
changes of variables as in section 2, we arrive at:
\begin{eqnarray*}
- ES_{\alpha} &=& \frac{1}{\alpha } \sum_{i=1}^m \beta_{i}
\int_{\{ |\delta A_{i}| z_{1} \leq - \delta \cdot \mu_{i} -
VaR_{\alpha } \}}  \ (|\delta A|z_{1} +
\delta \cdot \mu_{i} ) \   g_{i}( {\| z\|}^{2} )  dz \\
&=& \frac{1}{\alpha } \sum_{i=1}^m \beta_{i} \Big[ \int_{\{
|\delta A| z_{1} \leq - \delta \cdot \mu_{i} - VaR_{\alpha } \}} \
|\delta A_{i}|z_{1} \ g_{i}( {\| z\|}^{2} ) \ dz \ + \ \delta
\cdot \mu_{i}\Big] .
\end{eqnarray*}
The final integral on the right hand side can be treated as
before, by writing $ {\| z\|}^{2} = z_{1}^2 + {\| z^{'}\|}^{2} $
and introducing spherical coordinates $z^{'} = r\xi $, $\xi\in
S_{n-2} $, leading to:
$$
- ES_{\alpha } = \sum_{i=1}^m \beta_{i} \delta \cdot \mu_{i}  +
\frac{|S_{n-2}|}{\alpha } \sum_{i=1}^m \beta_{i} \int_{0}^{\infty}
r^{n-2} \Big{[} \int_{ -\infty }^{\frac{- \delta \mu_{i}^{t}
-VaR_{\alpha }}{|\delta A_{i}| }}    {  |\delta A_{i}| \; z_{1}} \
g_{i}( z_{1}^2 + r^{2} ) dz_{1}\Big{]} dr
$$
We now first change $z_1 $ into $-z_1 $, and then introduce $u =
z_1^2 + r^2 $, as before. If we recall that, by theorem
\ref{VaR-elliptic},
$$
q_{\alpha,i}^g = {\frac{ \delta \cdot \mu_{i} +VaR_{\alpha
}}{|\delta A_{i}| }}
$$
then, simply writing $q_{\alpha,i} $ for $q_{\alpha,n}^{f_{X}} $,
we arrive at:
\begin{eqnarray*}
ES_{\alpha } &=& - \sum_{i=1}^m \beta_{i} \Big[ \delta \cdot
\mu_{i} + |\delta A | \; \frac{|S_{n - 2 } | }{\alpha } \cdot \int
_{q_{\alpha,i } }^{\infty} \int _{z_1 ^2 } ^{\infty } z_1 (u - z_1
^2 )
^{\frac{n - 3 }{2 } } g(u) \ du \ dz_1 \Big] \\
&=& - \sum_{i=1}^m \beta_{i} (\delta \cdot \mu) + \sum_{i=1}^n
\beta_{i} |\delta A_{i} | \; \frac{|S_{n - 2 } | }{\alpha } \cdot
\int _{q_{\alpha,i } ^2 } ^{\infty } \frac{1 }{n - 1 } \left( u -
q_{\alpha, i } ^2 \right) ^{\frac{n - 1 }{2 } } \ g_{i}(u) \ du ,
\end{eqnarray*}
since
$$
\int _{q_{\alpha,i } } ^{\sqrt{u } } z_1 \left( u - z_1 ^2 \right)
^{\frac{n - 3 }{2 } } dz_{1} = \frac{1 }{n - 1 } \left( u -
q_{\alpha,i } ^2 \right) ^{\frac{n - 1 }{2 } } .
$$
After substituting the formula for $|S_{n - 2 } | $ and using the
functional equation for the $\Gamma $-function, $\Gamma (x + 1 ) =
x \Gamma (x) $, we arrive at the following result:

\begin{theorem}
Suppose that the portfolio is linear in the risk-factors $X = (X_1
, \cdots , X_n )$: $\Delta \Pi=\delta \cdot X $ and that $X \sim
N(\mu ,\Sigma ,\phi ) $, with pdf $f_{X}(x)= \sum_{i=1}^n
\beta_{i} {|\Sigma_{i}|}^{-1/2} g((x-\mu_{i}
)\Sigma_{i}^{-1}(x-\mu_{i} )^{t}) $. If we replace $q_{\alpha } $
by his value, then the expected Shortfall at level $\alpha $ is
given by : \begin{equation} ES_{\alpha } = - \sum_{i=1}^m \beta_{i} ( \delta
\cdot \mu_{i}) + \sum_{i=1}^m \beta_{i} |\delta \Sigma_{i} \delta
^t |^{1/2 } \cdot \frac{\pi^{\frac{n-1}{2}}}{\alpha \cdot
\Gamma(\frac{n + 1}{2})} \cdot \int _{(q_{\alpha , i } ^g )^2 }
^{\infty } \left( u - (q_{\alpha , i } ^g )^2 \right) ^{\frac{n -
1 }{2 } } \ g_{i}(u) \ du . \label{EESformula} \end{equation}
\end{theorem}
\begin{remark}\rm{
If we are in situations where $\mu=\mu_i$ and $\Sigma_{i}=\Sigma$
for all $i=1,\ldots,n$, therefore $q_{\alpha,i}$ does not depend
to $i$. It will depend only to the $q_{\alpha}$ given by the
mixture of elliptic VaR. In effect,
$q_{\alpha,i}=q_{\alpha}=q_{\alpha }^{ME-VaR_{\alpha}}$ such that
$$ VaR _{\alpha } = - \delta \cdot \mu +
q_{\alpha }^{ME-VaR_{\alpha}} \cdot \sqrt{ \delta \Sigma \delta ^t
}.$$}
\end{remark}
We therefore obtain the following corollary:
\begin{corollary}
Suppose that the portfolio is linear in the risk-factors $X = (X_1
, \cdots , X_n )$: $\Delta \Pi=\delta \cdot X $ and that $X \sim
N(\mu ,\Sigma ,\phi ) $, with pdf $f_{X}(x)= \sum_{i=1}^m
\beta_{i} {|\Sigma|}^{-1/2} g_{i}((x-\mu )\Sigma^{-1}(x-\mu )^{t})
$. If we replace $q_{\alpha } $ by his value, then the expected
Shortfall at level $\alpha $ is given by :
 \begin{equation}
 ES_{\alpha } = -  \delta \cdot \mu   \ + \ q_{\alpha}^{ME-ES} \cdot \sqrt{\delta \Sigma \delta ^t }
\end{equation} where \begin{equation} q_{\alpha}^{ME-ES}=\frac{\pi^{\frac{n-1}{2}}}{\alpha
\cdot \Gamma(\frac{n + 1}{2})}\sum_{i=1}^m \beta_{i} \cdot \int _{
(q_{\alpha }^{ME-VaR}) ^2 } ^{\infty } \left( u - (q_{\alpha
}^{ME-VaR}) ^2 \right) ^{\frac{n - 1 }{2 } } \ g_{i}(u) \ du .
\label{ESELL} \end{equation}
\end{corollary}

\subsection{Application: Mixture of Student-$\it{t}$ Expected Shortfall}

\medskip
In the case of multi-variate t-student distributions we have that
$g_{i}(u)= C(\nu_{i} ,n) {(1+u/\nu_{i} )}^{-\frac{(n+\nu_{i}
)}{2}} $, with $C(\nu_{i} , n ) $ given in section 2. Let us
momentarily write $q $ for $q_{\alpha , \nu_{i} }^{\mathfrak{t } } $.
Following \cite{S2}, we can evaluate the integral as follows:
\begin{eqnarray*}
&\ & \int _{q^2 } ^{\infty } (u - q)^{\frac{n - 1 }{2 } } \left( 1
+ \frac{u
}{\nu_{i} } \right) ^{-\frac{n + \nu_{i} }{2 } } du \\
&=& \nu_{i} ^{\frac{n + \nu_{i}}{2 } } (q^2 + \nu_{i}
)^{-(\frac{\nu_{i} - 1 }{2 } ) } B\left( \frac{\nu_{i} - 1 }{2 } ,
\frac{n + 1 }{ 2 } \right) .
\end{eqnarray*}
If we pose that : $$q_{\alpha}^{MST-ES}= \frac{1 }{\alpha \cdot
\sqrt{\pi } } \sum_{i=1}^m \beta_{i} \frac{\Gamma \left(
\frac{\nu_{i} - 1 }{2 } \right) }{\Gamma \left( \frac{\nu_{i} }{2
} \right) } \nu_{i} ^{\nu_{i} / 2 } \left( (q_{\alpha , \nu_{i} }
^{\mathfrak{t } } )^2 + \nu_{i} \right) ^{ \frac{1-\nu_{i} }{2 }
}\label{es}
$$
 After substitution in (\ref{EESformula}), we find, after some computations, the
following result:

\begin{theorem} The Expected Shortfall at confidence level $1 -
\alpha $ for a multi-variate Student-distributed linear portfolio
$\delta \cdot X $, with
$$
f_{X}(x)= \sum_{i=1}^m \beta_{i} \frac{\Gamma (\frac{\nu_{i} +
n}{2})}{\Gamma(\nu_{i}/2).
 \sqrt{|\Sigma|(\nu_{i} \pi)^n }}  {\Big{(}1+\frac{(x-\mu)^{t}
 \Sigma^{-1}(x-\mu) }{\nu_{i}} \Big{)}}^{-(\frac{\nu_{i} + n}{2})}  ,
$$
is given by:
\begin{eqnarray*}
ES_{\alpha , \nu } ^{\mathfrak{t } } &=& - \delta \cdot \mu
\sum_{i=1}^m \beta_{i} + |\delta \Sigma \delta ^t |^{1/2 }
\sum_{i=1}^m \beta_{i} \frac{1 }{\alpha \cdot \sqrt{\pi } }
\frac{\Gamma \left( \frac{\nu_{i} - 1 }{2 } \right) }{\Gamma
\left( \frac{\nu_{i} }{2 } \right) } \nu_{i} ^{\nu_{i} / 2 }
\left( (q_{\alpha }^{MST-VaR})^2 + \nu_{i} \right) ^{\left(
\frac{1-\nu_{i}}{2 }
\right) } \\
&=& -   \delta \cdot \mu  +   |\delta \Sigma \delta ^t |^{1/2 }
\sum_{i=1}^{m} \frac{\beta_{i} }{\alpha \cdot \sqrt{\pi } }
\frac{\Gamma \left( \frac{\nu_{i} - 1 }{2 } \right)\nu_{i}
^{\nu_{i} / 2 } }{\Gamma \left( \frac{\nu_{i} }{2 } \right) }
\left( \left( {\frac{ \delta \cdot \mu + VaR_{\alpha }}{|\delta
\Sigma \delta |^{1/2 } }} \right)^2 + \nu_{i} \right)
^{\left( \frac{1-\nu_{i}  }{2 } \right) } \\
&=& - \delta \cdot \mu  \ + \  q_{\alpha}^{MST-ES}\cdot
\sqrt{\delta \Sigma \delta ^t }
\end{eqnarray*}

\end{theorem}
\noindent The Expected Shortfall for a linear Student portfolio is
therefore given by a completely explicit formula, once the VaR is
known. Observe that, as for the VaR, the only dependence on the
portfolio dimension is through the portfolio mean $\delta \cdot
\mu $ and the portfolio variance $\delta \Sigma \delta ^t . $
\medskip
\subsection{Some Numerical Result of Delta Mixture-Student Expected Shortfall}
\medskip
\medskip

Here we give some numerical results when applying the corollary
\ref{CorrmixtVaR2}, in the situation where $m=2$.

For given $s=q_{\beta,\nu_{1},\nu_{2}}=q_{\alpha}^{MS-VaR}$, which
is the solution of
$$F(s,\beta,\nu_{1},\nu_{2})=\alpha,\label{mixtwostu2}$$
by introducing the fonction H such that
 \begin{equation}
H(s,\beta,\nu_{1},\nu_{2})= \beta \cdot  H_{1}(s) + (1-\beta)\cdot
H_{2}(s), \end{equation} where

$$H_{j}(s)=\frac{\beta_{i} }{\alpha \cdot \sqrt{\pi } } \frac{\Gamma
\left( \frac{\nu_{i} - 1 }{2 } \right) }{\Gamma \left(
\frac{\nu_{i} }{2 } \right) } \nu_{i} ^{\nu_{i} / 2 } \left( s^2 +
\nu_{i} \right) ^{ \frac{1-\nu_{i} }{2 }}$$

  for $j=1,2$. For given  as inputs
$\beta$, $\nu_{1}$ and $\nu_{2}$, we  give a table that contains
some values of
$q_{\alpha}^{MST-ES}=H(q_{\alpha}^{MST-VaR},\beta,\nu_{1},\nu_{2})=q_{\beta,\nu_{1},\nu_{2}}^{MST-ES}$.
\medskip
\begin{enumerate}
\item In the case where $\alpha=0.01$, we obtain some solutions of
(\ref{mixtwostu2}) in the following table:
\begin{center}
\begin{tabular}{|c|c|c|c|c|c|}

  \hline

$(\nu_{1},\nu_{2})$  & (2,3)  & (3,4) & (4,6) & (7,15) & (8,40)  \\
\hline
   $q_{0.25,\nu_{1},\nu_{2}}^{MST-ES}$   & 6.36587 & 1.29375  & 0.243125 & 0.00290856 & 0.000681262   \\
 \hline
 $q_{0.30,\nu_{1},\nu_{2}}^{MST-ES}$   & 7.01881 & 1.41000  & 0.279435 & 0.00341273 & 0.000793844  \\
 \hline
$q_{0.35,\nu_{1},\nu_{2}}^{MST-ES}$   & 7.64714 & 1.52252  & 0.31424 & 0.00389277 & 0.0008997532  \\
 \hline
 $q_{0.40,\nu_{1},\nu_{2}}^{MST-ES}$   & 8.25196 & 1.63141  & 0.34759 & 0.0043495 & 0.000997532  \\
 \hline
 $q_{0.45,\nu_{1},\nu_{2}}^{MST-ES}$   & 8.83444 & 1.73679  & 0.379538 &  0.00478369  & 0.00108926    \\
 \hline
 $q_{0.50,\nu_{1},\nu_{2}}^{MST-ES}$   & 9.3957 & 1.83877  & 0.410131 & 0.00519619 & 0.00117468    \\
 \hline
\end{tabular}
\end{center}
\medskip

\item In the case where $\alpha=0.001$, we obtain some solutions
of (\ref{mixtwostu2}) in the following table:
\begin{center}
\begin{tabular}{|c|c|c|c|c|c|}

  \hline

$(\nu_{1},\nu_{2})$  & (2,3)  & (3,4) & (4,6) & (7,15) & (8,40)  \\
\hline
   $q_{0.25,\nu_{1},\nu_{2}}^{MST-ES}$   & 20.8961 & 3.03289  & 0.576689 & 0.00661826 & 0.00164597   \\
 \hline
 $q_{0.30,\nu_{1},\nu_{2}}^{MST-ES}$   & 23.1642 & 3.32289  & 0.666054 & 0.0074621 & 0.00180969  \\
 \hline
$q_{0.35,\nu_{1},\nu_{2}}^{MST-ES}$   & 25.2707 & 3.58757  & 0.716427 & 0.008196 & 0.00194229  \\
 \hline
 $q_{0.40,\nu_{1},\nu_{2}}^{MST-ES}$   & 27.239 & 3.83719  & 0.776394 & 0.00883632 & 0.00205071  \\
 \hline
 $q_{0.45,\nu_{1},\nu_{2}}^{MST-ES}$   & 29.0885 & 4.07077  & 0.830853 &  0.00939711  & 0.00214048    \\
 \hline
 $q_{0.50,\nu_{1},\nu_{2}}^{MST-ES}$   & 30.8351 & 4.28993  & 0.880508 & 0.00989055 & 0.00221577    \\
 \hline
\end{tabular}
\end{center}
\medskip

\end{enumerate}

\subsection{Delta-Theta Approximation of a Portfolio}

\medskip
In the case where we dealt with portfolio that contains
derivatives, we will consider the Greek Theta of the portfolio by
replace the Delta approximation known in financial literature by
the {\em{ Delta-Theta approximation}}.

In clear, suppose that we are holding a portfolio of derivatives
depending on $n $ underlying assets $X(t)=(X_{1}(t),\ldots,
X_{n}(t))$, with elliptically distributed log-returns $r_j
=log(X_{j}(t)/X_{j}(0)) $, over some fixed small time-window
[0,t]. The portfolio's present value $V $ will in general be some
complicated non-linear function of the $X_{i} $'s. To obtain a
first approximation of its VaR, we simply approximate the present
Value V of the position using a first order Taylor expansion:
\begin{eqnarray}
V(X(t),t)&\approx& V(X(0),0) + \sum_{i=1}^n  \frac{\partial V}{
\partial X_{i}} (X(0),0) (X_{i}(t)-X_{i}(0))+ t\cdot  \frac{\partial V}{
\partial{t}}(X(0),0)\nonumber \\
&=& V(X(0),0)+ \sum_{i=1}^n  \frac{\partial V}{
\partial X_{i}} (X(0),0)  X_{i}(0) \left( \exp( r_i) -1 \right) +
\Theta \cdot t \nonumber \\
&\approx& V(X(0),0)+ \sum_{i=1}^n  \delta_{i}  r_i + \Theta \cdot
t \label{Delta-Theta}
\end{eqnarray}

 From this, we can then approximate the Profit \& Loss function
as
$$
\Delta V \approx \delta \cdot r^t + \Theta\cdot t ,
$$
where we put $r=(r_{1},\ldots,r_{n})$ and
$\delta=(\delta_{1},...,\delta_{n})$ with
$\delta_{i}=X_{i}(0)\cdot\frac{\partial V}{ \partial
X_{i}}(X(0),0)$. The entries of the $\delta$ vector are called the
"delta equivalents " for the position, and they can be interpreted
as the sensitivities of the position with respect to changes in
each of the risk factors.
   In this particular case, we have substitute the {{\em{Delta normal VaR}} as known in the financial literature, by
   the {{\em{Delta-Theta Elliptic VaR}}  given by the following corollary of the theorem (\ref{VaR-elliptic}) :

\begin{corollary} \label{Delta-Theta-elliptic}
Suppose that the portfolio's Profit \& Loss function over the time
window of interest is, to good approximation, given by $\Delta
\Pi=\delta \cdot r^t + \Theta \cdot t$ , with constant portfolio
weights $\delta=(\delta _1 , \ldots,\delta_{n} ) $. Suppose
moreover that the random vector $r = (r_1 , \cdots , r_n ) $ of
underlying log-returns follows a continuous elliptic distribution,
with probability density given by $f_r (x)= {|\Sigma|}^{-1/2}
g((x-\mu )\Sigma^{-1}(x-\mu )^{t}) $ where $\mu$ is the vector
mean and $\Sigma $ is the variance-covariance matrix, and where we
suppose that $g (s^2 ) $ is integrable over $\mathbb{R } $,
continuous and nowhere 0. Then the portfolio's {{\em
Delta-Theta-elliptic VaR}} $VaR _{\alpha } $ at confidence $1 -
\alpha $ is given by
$$
VaR _{\alpha } =  - \delta \cdot \mu^t + \Theta \cdot t +
q_{\alpha,n } ^g \cdot \sqrt{ \delta \Sigma \delta ^t } ,
$$
where $s = q_{\alpha,n } ^g $ is the unique positive solution of
the transcendental equation
$$
\alpha = G ( q_{\alpha,n }^g ) .
$$
\end{corollary}

The Expected Shortfall of such portfolios is given by the
following corollary
\begin{corollary} \label{Delta-Theta-EllipticES}
Suppose  that the portfolio's Profit \& Loss function over the
time window of interest is, to good approximation, given by
$\Delta \Pi=\delta \cdot r^t + \Theta \cdot t$ ,and that $r \sim
N(\mu ,\Sigma ,\phi ) $, with pdf $f(x)= {|\Sigma|}^{-1/2}
g((x-\mu )\Sigma^{-1}(x-\mu )^{t}) $, then the {\em Delta-Theta
Elliptic Expected Shortfall } or {\em Delta-Theta ES} at
confidence level $1-\alpha $ is  :
 \begin{equation}
ES_{\alpha } = - \delta \cdot \mu^t + \Theta \cdot t + |\delta
\Sigma \delta ^t |^{1/2 } \cdot \frac{\pi^{\frac{n-1}{2}}}{\alpha
\cdot \Gamma(\frac{n + 1}{2})} \cdot \int _{(q_{\alpha , n } ^g
)^2 } ^{\infty } \left( u - (q_{\alpha , n } ^g )^2 \right)
^{\frac{n - 1 }{2 } } \ g(u) \ du .\end{equation}
\end{corollary}

\begin{remark} \rm{In short-term Risk Management, one can usually
assume that $\mu \simeq 0 $. In that case, for $t=1$ we have
$$VaR_{\alpha}=\Theta+ \sqrt{\delta \Sigma \delta^{t}} \cdot q_{\alpha,n } ^g ,$$
$$ES_{\alpha } =  \Theta  + |\delta \Sigma \delta ^t
|^{1/2 } \cdot \frac{\pi^{\frac{n-1}{2}}}{\alpha \cdot
\Gamma(\frac{n + 1}{2})} \cdot \int _{(q_{\alpha , n } ^g )^2 }
^{\infty } \left( u - (q_{\alpha , n } ^g )^2 \right) ^{\frac{n -
1 }{2 } } \ g(u) \ du .$$
 }
\end{remark}
\noindent As before, The preceding can immediately be specialized
to a Student $\mathfrak{t } $-distributions to estimate the {\em
Delta-Theta Student VaR} and the {\em Delta-Theta Student ES}. The
details will be left to the reader.
\subsection{Portfolios of Equities}
\medskip

A special case of the preceding is that of an equity portfolio,
build of stock $S_1 , \ldots , S_n $ with  joint log-returns
$r=(r_1 (t) , \ldots , r_n (t)) . $ In this case, the portfolio's
Profit \& Loss function over the time window [0,t] of interest is,
to good approximation, given by
\begin{eqnarray*}
\Pi(t) - \Pi (0) &=& \sum_{i=1}^n w_{i} S_{i}(0) (
S_{i}(t)/S_{i}(0) -1 ) \\
&\approx & \sum _{i = 1 } ^n w_i S_i (0) r_i (t)=\delta \cdot r^t,
\end{eqnarray*}
where this approximation will be good if the $r_i (t) $ are small.
In this case the preceded  theorems are applicable where
$\delta=(w_1 S_1 (0),\ldots,w_n S_n (0))$ and
$r_{j}(t)=log(X_{j}(t)/X_{j}(0)) $ for j=1,\ldots,n.

\subsection{Businesses as Linear Portfolios of Business Units} An
interesting way of looking upon an big enterprize, e.g. a
multi-national or a big financial institution, is by considering
it as a sum of its individual business units, cf. Dowd \cite{DO}.
If $X_{j}$, is the variation of price or of profitability of
business unit j in one period, then the variation of price of the
agglomerate in the same period will be
$$
\Delta \Pi= X_{1}+\cdots+ X_{n} .
$$
The entire institution is therefore modelled by a linear
portfolio, with $\delta=(1, 1, \ldots , 1) $, to which the results
of this paper can be applied, if we model the vector of individual
price variations by a multi-variate elliptic distribution. VaR,
incremental VaR (see below) and Expected Shortfall will be
relevant here. For more details see Dowd \cite{DO}, chapter XI .

\subsection{Incremental VaR}
Incremental VaR is defined in \cite{MX}  as the statistic that
provides information regarding the  sensitivity of VaR to changes
in the portfolio holdings. It therefore gives an estimation of the
change in VaR resulting from a risk management decision. Results
from \cite{MX} for incremental VaR with normally distributed
risk-factors generalize straightforwardly to elliptically
distributed ones: if we denote by $IVaR_{i}$ the incremental VaR
for each position in the portfolio, with $\theta_{i}$ the
percentage change in size of each position, then the change in VaR
will be given by
$$
\Delta VaR=\sum \theta_{i}IVaR_{i} \label{ivar1}
$$
By using the definition of $IVaR_{i}$ as in \cite{MX}  (2001), we
have that \begin{equation} IVaR_{i}=\omega_{i}\frac{\partial VaR }{\partial
\omega_{i}} \end{equation} with $\omega_{i}$ is the amount of money invested
in instrument $i $. In the case of an equity portfolio in the
elliptically distributed assets, we have seen that, assuming
$\mu=0$,
$$
VaR_{\alpha} = - q_{\alpha,n}^g \sqrt{\delta \Sigma \delta^t } ,
$$
We can then calculate $IVaR_{i}$ for the i-th constituent of
portfolio as
$$
IVaR_{i}=\omega_{i}\frac{\partial VaR }{\partial
\omega_{i}}=\omega_{i}\gamma_{i}
$$
with
$$
\gamma= - q_{\alpha,n}^g \frac{\Sigma \omega}{\sqrt{\delta \Sigma
\delta^t } } .
$$
The vector $\gamma $ can be interpreted as a gradient of
sensitivities of VaR with respect to the risk factors. This is the
same as in \cite{MX}, except of course that the quantile has
changed from the normal one to the one associated to $g . $

\subsection{Problem of the aggregation of risks}
Suppose that we have a constituted portfolio with several under
portfolios of assets from different markets.
  Given the Value-at-Risk of the portfolios constituting the global portfolio, under
  the hypothesis that the joined risks factors follow an {\em  elliptic
  distribution }, the question is how to get the VaR of the global
  portfolio.

  In order to be clearer and simpler, let us consider  a global constituted portfolio of 2
  under portfolios from different markets with respective weights $\mathbf{\delta_1} $ and $\mathbf{\delta_2} $.
   $\mathbf{\Sigma_1} $ represents  the matrix of interrelationship in the under portfolio of market 1;
    $\mathbf{\Sigma_2} $ represents
   the matrix of interrelationship in the  under portfolio of market 2.
  One will be able to write the matrix of interrelationship of a global portfolio like this:
  \begin{displaymath}
  \mathbf{\Sigma}=
  \left ( \begin{array}{cc}
  \mathbf{\Sigma_1} & \mathbf{\Sigma_{12}}\\
  \mathbf{\Sigma_{12}}^{t} & \mathbf{\Sigma_2 }  \\
  \end{array} \right),
  \end{displaymath}
 where $\mathbf{\Sigma_{12}}$ is the correlation matrix that takes
 into consideration the interaction between the market $ \bf{M_1} $ and
  the market $\bf{M_2}$ .   If $\mathbf{\delta}^t=( \mathbf{\delta_1},\mathbf{\delta_2}) $, we
  have
        \begin{equation} \mathbf{\delta}^t \mathbf{\Sigma}\mathbf{\delta}=
        \mathbf{\delta_1}^t \mathbf{\Sigma_1}\mathbf{\delta_1}+ \mathbf{\delta_2}^t
        \mathbf{\Sigma_2}\mathbf{\delta_2}+ 2 \cdot \mathbf{\delta_1}^t
        \mathbf{\Sigma_{12}}\mathbf{\delta_2} .\label{agreg1}
        \end{equation}
 Therefore, since we know that when  $\mu \approx 0 $,
    we have
    $$
VaR _{\alpha } =  q_{\alpha,n } ^g \cdot \sqrt{ \delta \Sigma
\delta ^t } ,
$$
 the Value-at-Risk of the  global  portfolio will be given
by \begin{equation} VaR _{\alpha }(\bf{M}) = \sqrt{ VaR _{\alpha }(\bf{M_1})^2
+ VaR _{\alpha }(\bf{M_2})^2 + 2 [ q_{\alpha,n } ^g ]^2 \cdot
\mathbf{\delta_1}^t
        \mathbf{\Sigma_{12}}\mathbf{\delta_2}}.
   \label{agreg2}
     \end{equation}

     An implicit interrelationship with the hypothesis of {\em elliptic distribution} is
     obtained in an analogous way, like in the case  where one
     works with the hypothesis of the normal distribution. Note that, one will distinguish several
        situations from  the behavior of $\mathbf{\Sigma_{12}}$.
      With  some simple operations, the implicit
      interrelationship is
      \begin{equation}
      \phi=\frac{\mathbf{\delta_1}^t
        \mathbf{\Sigma_{12}}\mathbf{\delta_2}}{\sqrt{(\mathbf{\delta_1}^t \mathbf{\Sigma_1}\mathbf{\delta_1})(\mathbf{\delta_2}^t
        \mathbf{\Sigma_2}\mathbf{\delta_2})}}
        \end{equation}
        with the Value-at-Risk $VaR _{\alpha }(\bf{M})$ of the global portfolio being given as follows:
\begin{equation} VaR _{\alpha }(\bf{M}) = \sqrt{ [VaR _{\alpha }(\bf{M_1})]^2 +
[VaR _{\alpha }(\bf{M_2})]^2 + 2 \phi \cdot VaR_{\alpha}(\bf{M_1})
VaR _{\alpha }(\bf{M_2}))}.
   \end{equation}

   Also, for $\mu \approx 0 $, $$ES_{\alpha } = K_{ES,\alpha}^g
\cdot \sqrt{\delta \Sigma \delta ^t} ,$$ therefore by using the
same technics that proves (\ref{agreg2}), we have that the
expected shortfall of the global portfolio is given by:

 \begin{equation} ES
_{\alpha }(\bf{M}) = \sqrt{ ES _{\alpha }(\bf{M_1})^2 + ES
_{\alpha }(\bf{M_2})^2 + 2 [ K_{ES,\alpha}^g ]^2 \cdot
\mathbf{\delta_1}^t
        \mathbf{\Sigma_{12}}\mathbf{\delta_2}}.
   \label{agregES1}
     \end{equation}
This imply that

 \begin{equation} ES _{\alpha }(\bf{M}) = \sqrt{ [ES _{\alpha
}(\bf{M_1})]^2 + [ES _{\alpha }(\bf{M_2})]^2 + 2 \phi_{ES} \cdot
ES_{\alpha}(\bf{M_1}) ES _{\alpha }(\bf{M_2}))},
  \label{agregES2}
   \end{equation}
where
 \begin{equation}
      \phi_{ES}=\frac{\mathbf{\delta_1}^t
        \mathbf{\Sigma_{12}}\mathbf{\delta_2}}{\sqrt{(\mathbf{\delta_1}^t \mathbf{\Sigma_1}\mathbf{\delta_1})(\mathbf{\delta_2}^t
        \mathbf{\Sigma_2}\mathbf{\delta_2})}}
        \label{agregES3}
        \end{equation}
\begin{remark}{\rm
     The result about the agregation of risks work so well in the
     situation where, the joint risk factors of our portfolio
     changes with mixture of ellitic distributions as define in
     (\ref{dens}), and where all $\Sigma_{i}=\Sigma$, for
    $ i=1,\ldots,m$. In particular, when $\mu_{i}=\mu$, we have the results (\ref{agregES2}) and (\ref{agreg2})}.
\end{remark}

\section{conclusion}
\medskip
In this paper, we have shown how to reduce the estimation of
Value-at-Risk for linear elliptic portfolios to the evaluation of
one dimensional integrals which, for the special case of a mixture
of  $t$-Student distributions, can be explicitly evaluated in
terms of a hypergeometric function. We have also given a similar,
but simpler, integral formula for the expected shortfall of such
portfolios which, again, can be completely evaluated in the
Student case. Following the calculations in the case of Delta
mixture-Student VaR, we indicated how to extend it to the case of
mixture of t-distributions expected shortfall . We finally
surveyed some potential application areas.

\end{document}